\documentclass[a4paper,12pt]{amsart}

\usepackage{amssymb,amscd,amsthm,mathrsfs}
\usepackage[ansinew]{inputenc}
\usepackage[centertags]{amsmath}

\newcommand{\en}{\subset}
\newcommand{\N}{\mbox{$\mathbb{N}$}}

\newcommand{\Z}{\mbox{$\mathbb{Z}$}}

\newcommand{\R}{\mbox{$\mathbb{R}$}}
\newcommand{\Res}{\mbox{$\mathcal{R}$}}
\newcommand{\T}{\mbox{$\mathbb{T}$}}
\newcommand{\cO}{\mathcal{O}}
\newcommand{\cC}{\mathcal{C}}
\newtheorem{teo}{Theorem}[section]
\newtheorem*{teo*}{Theorem}
\newtheorem*{teoA1}{Theorem A.1}
\newtheorem*{teoA2}{Theorem A.2}

\newtheorem*{teoB1}{Theorem B.1}
\newtheorem*{lemaB1}{Lemma B.1}
\newtheorem*{lemaB2}{Lemma B.2}
\newtheorem*{cor*}{Corollary}
\newtheorem{conj}{Conjecture}
\newtheorem{quest}{Question}
\newtheorem{claim}[]{Claim}
\newtheorem*{claim*}{Claim}
\newtheorem{cor}{Corollary}[section]
\newtheorem{lema}{Lemma}[section]
\newtheorem{prop}{Proposition}[section]

\newcommand{\bi}{\begin{itemize}}
\newcommand{\ei}{\end{itemize}}

\theoremstyle{definition}

\theoremstyle{remark}
\newtheorem{obs}{Remark}[section]

\newcommand{\D}{\mbox{$\mathbb{D}$}}

\newcommand{\eps}{\varepsilon}

\newcommand{\dem}[1]{\vspace{.05in}{\sc\noindent Proof #1.}}
\newcommand{\lqqd}{\par\hfill {$\Box$} \vspace*{.1in}}
\newcommand{\finobs}{\par\hfill{$\diamondsuit$} \vspace*{.1in}}

\newcommand{\cA}{\mathcal A}
\newcommand{\cB}{\mathcal B}

\newcommand{\U}{\mathcal{U}}
\newcommand{\trans}{\mbox{$\,{ \top} \;\!\!\!\!\!\!\raisebox{-.3ex}{$\cap$}\,$}}
\newcommand{\umb}[1]{\hat #1}
\DeclareMathOperator{\Diff}{Diff}
\DeclareMathOperator{\Tang}{Tang}
\DeclareMathOperator{\diametro}{diam}

\newcommand{\diff}[1]{\Diff^{#1}(M)}

\author[Rafael Potrie]{Rafael Potrie}
\address{CMAT, Facultad de Ciencias, Universidad de la Rep\'ublica, Uruguay}\address{LAGA; Institut Galilee, Universite Paris 13, Villetaneuse, France}
\email{rpotrie@cmat.edu.uy}

\title{Generic bi-Lyapunov stable homoclinic classes}
\thanks{The autor was partially supported by ANR Blanc DynNonHyp BLAN08-2$\_$313375 and ANII, FCE2007 n577}

\begin{document}



\begin{abstract}
We study, for $C^1$ generic diffeomorphisms, homoclinic classes which are Lyapunov stable both for backward and forward iterations. We prove that they must admit a dominated splitting and show that under some hypothesis they must be the whole manifold. As a consequence of our results we also prove that in dimension 2 the class must be the whole manifold and in dimension 3, these classes must have nonempty interior. Many results on Lyapunov stable homoclinic classes for $C^1$-generic diffeomorphisms are also deduced. %
\end{abstract}

\maketitle

\section{Introduction}

The study of the dynamics of chain recurrence classes and their interaction with other chain recurrence classes has become a major problem in $C^1$-generic dynamics since \cite{BC}, where generic dynamics and Conley's theory (see section 10.1 of \cite{Rob}) were unified. In particular, this has raised interest in the study of homoclinic classes, which are generically the chain recurrence classes containing periodic orbits.

So, an interesting problem is to study the possible structures and dynamics a homoclinic class may have, particularly those related to the global dynamics of the diffeomorphism. At the moment, very little is known related to this problem, specially when the class is \emph{wild}, that is, accumulated by infinitely many distinct chain recurrence classes (isolated classes are now quite well understood, see \cite{beyondhip} Chapter 10). In particular, very natural and simple questions remain wide open such as if one of these homoclinic classes may have nonempty interior and not be the whole manifold (the progress made so far has to do with \cite{ABCD},\cite{ABD},\cite{PotS}).

On the other hand, also a lot of work has been done towards the understanding of generic dynamics far from homoclinic tangencies (see for example \cite{PujSam},\cite{W}, \cite{C},\cite{Y}), but also very few is known in their presence.

We deal in this paper with this kind of difficulties in a simpler context than just wild homoclinic classes. We study homoclinic classes which are \emph{bi-Lyapunov stable}, in particular, which are saturated by stable and unstable sets (for example, homoclinic classes with non empty interior have this property). We are able, by using new techniques developed in \cite{Go3}, to prove that these classes admit a dominated splitting, a weak form of hyperbolicity which helps to order the dynamics in a class (see \cite{beyondhip}, Appendix B). Also, we handle the case where the class is far from tangencies and we prove that in that case, if a generic diffeomorphism admits a homoclinic class which is bi-Lyapunov stable then the diffeomorphism must be transitive. Our results respond affirmatively to the second part of Problem 5.1 in \cite{ABD} and to the other part also in dimension 2 and 3.

The main technical novelty of this paper relies in the use of Lyapunov stability together with a recent result from \cite{Go3} to control the bifurcation of periodic points.

Finally, we would like to add that many of our results treat the case of homoclinic classes which are only Lyapunov stable for $f$ (and not for $f^{-1}$) and could have independent interest for other applications.

\subsection{Context of the problems}

\smallskip

Let $M$ be a compact connected boundaryless manifold of dimension
$d$ and let $\diff 1$ be the set of diffeomorphisms of $M$ endowed
with the $C^1$ topology.

We shall say that a property is $C^1$-\emph{generic}, or just \emph{generic}, if
and only if there exists a residual set ($G_{\delta}-$ dense) $\Res$ of $\diff 1$ such that every $f \in \Res$ satisfies that property. Some well known results in the theory of generic dynamical systems are presented in the Appendix A and they will be referred to when they are used.

For a hyperbolic periodic point $p \in M$ of some diffeomorphism $f$
we denote its \emph{homoclinic class} by $H(p,f)$. It is defined as the
closure of the transversal intersections between the stable and
unstable manifolds of the orbit of $p$. Homoclinic classes can also be defined in terms of homoclinically related periodic points by defining $H(p,f)$ as the closure of the set of the homoclinically related periodic orbits, both definitions coincide. We recall that two hyperbolic periodic orbits $\cO,\cO'$ are said to be \emph{homoclinically related} if $W^s(\cO)\trans W^u(\cO') \neq \emptyset$ and $W^s(\cO')\trans W^u(\cO)\neq \emptyset$.

It is a very important problem to study the structure of homoclinic classes since they are basic pieces of the dynamics (see \cite{beyondhip}, Chapter 10). In particular, a very natural question is if they can admit non-empty interior. The only known examples are robustly transitive diffeomorphisms for which there is (generically) only one homoclinic class which coincides with the whole manifold (\cite{beyondhip}, Chapter 7).

In \cite{ABD} the following conjecture was posed (it also appeared as Problem 1 in \cite{BC})

\begin{conj}[\cite{ABD}]\label{Conjetura} There exists a residual set $\Res$ of $\diff 1$ of diffeomorphisms such that if $f\in \Res$ admits a homoclinic class with nonempty interior, then the diffeomorphism is transitive.
\end{conj}

Some progress has been made towards the proof of this conjecture
(see \cite{ABD},\cite{ABCD} and \cite{PotS}), in particular, it has been proved
in \cite{ABD} that isolated homoclinic classes as well as homoclinic
classes admitting a strong partially hyperbolic splitting
 verify the conjecture. Also, they proved that a homoclinic class with non empty
interior must admit a dominated splitting (see Theorem 8 in
\cite{ABD}).

In \cite{ABCD} the conjecture was proved for surface
diffeomorphisms, other proof for surfaces (which does not use the approximation by $C^2$ diffeomorphisms) can be found in \cite{PotS} where the codimension one case is studied.

Also, from the work of Yang (\cite{Y}) one can deduce the conjecture in the case $f$ is $C^1$-generic and far from homoclinic tangencies (we shall extend this remark in section \ref{sectionY}).

When studying some facts about this conjecture, in \cite{ABD} it was proved that if a homoclinic class of a $C^1$-generic diffeomorphism has nonempty interior then this class should be bi-Lyapunov stable. In fact, in \cite{ABD} they proved that isolated and strongly partially hyperbolic bi-Lyapunov stable homoclinic classes for generic diffeomorphisms are the whole manifold.

This concept is a priori weaker than having nonempty interior and it is natural to ask the following question.

\begin{quest}[Problem 1 of \cite{BC}]\label{Question} Is a bi-Lyapunov stable homoclinic class of a generic diffeomorphism necessarily the whole manifold?
\end{quest}

For generic diffeomorphisms, there exists a more general notion of basic pieces of the dynamics, namely, the chain recurrence classes (\footnote{The chain recurrent set is
 the set of points $x$ satisfying that for every $\eps>0$
 there exist an $\eps$-pseudo orbit form $x$ to $x$, that is,
 there exist points $x=x_0, x_1, \ldots x_k=x$ with $k\geq 1$ such that
 $d(f(x_i),x_{i+1})< \eps$. Inside the chain recurrent set, the chain recurrence classes are the equivalence classes of the relation given by $x\vdash \! \dashv y$ when for every $\eps>0$ there exists an $\eps-$pseudo orbit from $x$ to $y$ and one from $y$ to $x$ (see \cite{Rob}). A compact invariant set $K$ is chain transitive if $K$ is a chain recurrence class for $f_{/K}$.}). It is not difficult to deduce from \cite{BC} that, for generic diffeomorphisms, a chain recurrence class with non empty interior must be a homoclinic class, thus, the answer to Conjecture \ref{Conjetura} must be the same for chain recurrence classes and for homoclinic classes.

 However, we know that Question \ref{Question} admits a negative answer if posed for chain recurrence classes. Bonatti and Diaz constructed (see \cite{BD}) open sets of diffeomorphisms in every manifold of dimension $\geq 3$ admitting, for generic diffeomorphisms there, uncountably many bi-Lyapunov stable chain recurrence classes which in turn have no periodic points.

Although this may suggest a negative answer for Question \ref{Question} we present here some results suggesting an affirmative answer. In particular, we prove that the answer is affirmative for surface diffeomorphisms, and that in three dimensional manifold diffeomorphisms the answer must be the same as for Conjecture \ref{Conjetura}.

The main reason for which the techniques in \cite{ABCD} (or in \cite{PotS}) are not able to answer Question \ref{Question} for surfaces, is because differently from the case of homoclinic classes with interior, it is not so easy to prove that bi-Lyapunov stable classes admit a dominated splitting (in fact, the bi-Lyapunov stable chain recurrence classes constructed in \cite{BD} do not admit any). Here we are able, by using new techniques introduced by Gourmelon in \cite{Go3}, to show the existence of a dominated splitting for bi-Lyapunov stable homoclinic classes for generic diffeomorphisms. This will give us the result in dimension 2.

Also, using the same techniques, we are able to extend the results previously obtained in \cite{PotS} to the context of Question \ref{Question} which in turn allow us (combined with a new result of Yang, \cite{Y}) also to deduce an affirmative answer to the question in the far from tangencies context.


\subsection{Definitions and statement of results}

Let us first recall the definition of dominated splitting: a compact set
$H$, invariant under a diffeomorphism $f$, admits a \emph{dominated splitting}
if the tangent bundle over $H$ splits into two $Df$ invariant
subbundles, $T_HM=E\oplus F$, such that there exist $C>0$ and
$0<\lambda<1$, satisfying that for all $x\in H:$
$$\|Df^n_{/E(x)}\|\|Df^{-n}_{/F(f^n(x))}\|\le C\lambda^n.$$

In this case we say that the bundle $F$ dominates $E$. Let us remark that Gourmelon (\cite{GourmelonAdaptada}) proved that
there  always exists an adapted metric for which we can take $C=1$ in the definition. Given a dominated splitting $T_HM=E\oplus F$ we say that $\dim E$ is the \emph{index} of the dominated splitting. Similarly, for a periodic point $p$ of period $\pi(p)$, its \emph{index} is the dimension of the eigenspace associated to the set of eigenvalues of modulus smaller than $1$ for $Df_p^{\pi(p)}$. For a homoclinic class $H$, we define the \emph{minimal index} of $H$ to the minimum of the indexes of its periodic points (we define \emph{maximal index} analogously).

One can have dominated splittings into more than 2 subbundles (see \cite{beyondhip}, Appendix B), in particular, a splitting of the form $T_\Lambda M = \bigoplus_{i=1}^m E^i$ is dominated if for every $j<k$, $E^k$ dominates $E^j$.

We shall say that a bundle $E$ is \emph{uniformly contracting (expanding)} if there exists $n_0>0$ ($n_0<0$) such that $\|Df^{n_0}_{/E(x)}\|< 1/2$ $\forall x \in H$.

 Recall that a chain recurrent class $H$ is \emph{Lyapunov stable}\footnote{For $C^1$-generic diffeomorphisms, chain recurrence classes which are Lyapunov stable are also \emph{quasi-attractors}. See \cite{BC}.}\emph{ for $f$} if for every neighborhood $U$ of $H$ there is $V$ neighborhood of $H$ such that $f^n(\overline V) \en U$ $\forall n \geq 0$; in particular, it is easy to see that Lyapunov stability implies that $W^u(x) \en H$ for every $x\in H$. We shall say the class is \emph{bi-Lyapunov stable} if it is Lyapunov stable for $f$ and $f^{-1}$.

Given a hyperbolic periodic point $p$ in a homoclinic class $H$ of a diffeomorphism $f$ we have that its continuation is well defined in a small neighborhood $\U$ of $f$ and we denote $p_g$ to the continuation for $g\in \U$. Thus, we can also define the continuation of the homoclinic class $H$ for every $g \in \U$ as $H_g=H(g,p_g)$ the homoclinic class for $g$ of $p_g$. It is well known that given a hyperbolic periodic point $p$ of a diffeomorphism $f$ and $\U$ an open set containing $f$ where the continuation of $p$ is well defined, then, there exists a residual subset of $\U$ where the map $g \mapsto H(g,p_g)$ is continuous in the Hausdorff topology. Also, we can assume that if $H$ is Lyapunov stable (or bi-Lyapunov stable) for $f$, then $H(g,p_g)$ will also be Lyapunov stable (resp. bi-Lyapunov stable) for every generic $g$ in $\U$ (see Appendix A).

\begin{teo}\label{descdom} For every $f$ in a residual subset $\Res_1$ of $\diff 1$, if $H$ is a bi-Lyapunov stable homoclinic class for $f$, then, $H$ admits a dominated splitting. Moreover, it admits at least one dominated splitting with index equal to the index of some periodic point in the class.
\end{teo}

This theorem solves affirmatively the second part of Problem 5.1 in \cite{ABD}. We remark that Theorem \ref{descdom} does not imply that the class is not accumulated by sinks or sources. Also, we must remark that the theorem is optimal in the following sense, in \cite{BV} an example is constructed of a robustly transitive diffeomorphism (thus bi-Lyapunov stable) of $\T^4$ admitting only one dominated splitting (into two two-dimensional bundles) and with periodic points of all possible indexes for saddles.

The theorem is in fact stronger, we in fact prove the following theorem which may be interesting for itself dealing with Lyapunov stable homoclinic classes.

\begin{teo}\label{descdomLYAPUNOV} For every $f$ in a residual subset $\Res_1$ of $\diff 1$, if $H$ is a Lyapunov stable homoclinic class for $f$,  and there is a periodic point $p\in H$ of period $\pi(p)$ such that $det(Df^{\pi(p)}_p)\leq 1$, then, $H$ admits a dominated splitting.
\end{teo}

We recall now that a compact invariant set $H$ is \emph{strongly partially hyperbolic} if it admits a three ways dominated splitting $T_HM= E^s \oplus E^c \oplus E^u$, where $E^s$ is non trivial and uniformly contracting and $E^u$ is non trivial and uniformly expanding.

In the context of Question \ref{Question} it was shown in \cite{ABD} that generic bi-Lyapunov stable homoclinic classes admitting a strongly partially hyperbolic splitting must be the whole manifold. Thus, it is very important to study whether the extremal bundles of a dominated splitting must be uniform. We are able to prove this in the codimension one case, using only Lyapunov stability for future iterates.

\begin{teo}\label{maintheorem}
 Let $f$ be a diffeomorphism in a residual subset $\Res_2$ of $\diff 1$ with a homoclinic class $H$ which is Lyapunov stable
 admitting a codimension one dominated splitting $T_H M=E \oplus F$
 where $\dim(F)=1$. Then, the bundle $F$ is uniformly expanding for $f$.
\end{teo}

This theorem is an extension of a related result from \cite{PotS} (where Conjecture \ref{Conjetura} was studied) where the same result was proved in the case the class has nonempty interior.

As a consequence of this theorem we get the following easy corollaries.

\begin{cor}\label{Cor1} Let $H$ be a bi-Lyapunov stable homoclinic class
for a $C^1$-generic diffeomorphism $f$ such that $T_HM= E^1 \oplus E^2
\oplus E^3$ is a dominated splitting for $f$ and
$\dim(E^1)=\dim(E^3)=1$. Then, $H$ is strongly partially hyperbolic and $H=M$.
\end{cor}

\dem{\!\!} The class should be strongly partially hyperbolic because of
the previous theorem.  Corollary 1
of \cite{ABD} (page 185)  implies that $H=M$.
\lqqd

We say that a hyperbolic periodic point $p$ is \emph{far from tangencies} if there
is a neighborhood of $f$ such that there are no homoclinic tangencies associated to the stable and unstable manifolds of the continuation of $p$. The tangencies are of index $i$ if they are associated to a periodic point of index $i$, that is, its stable manifold has dimension $i$.
We  get the following result following \cite{ABCDW} (see also \cite{W} and \cite{Gou2}):

\begin{cor}\label{corolariolejostang} Let $H$ be a bi-Lyapunov stable homoclinic class
for a $C^1$-generic diffeomorphism $f$ which has a periodic point $p$ of index $1$ and a periodic point $q$ of index $d-1$ and such that $p$ and $q$ are far from tangencies . Then, $H=M$.
\end{cor}

\dem{\!\!} Using Corollary 3 of \cite{ABCDW} we are in the hypothesis of Corollary 1.1.
\lqqd

Using a result of Yang (Theorem 3 of \cite{Y}) and Theorem \ref{maintheorem} we are able to prove a similar result which is stronger than the previous corollary but which in turn, has hypothesis of a more global nature. We say that a diffeomorphism $f$ is \emph{far from tangencies} if it can not be approximated by diffeomorphisms having homoclinic tangencies for some hyperbolic periodic point. Notice that in the far from tangencies context, it is proved in \cite{Y} that a Lyapunov stable chain recurrence class must be an homoclinic class.

\begin{prop}\label{lejosdetangencias} There exists a $C^1$-residual subset of the open set of diffeomorphisms far from tangencies such that if $H$ is a bi-Lyapunov stable chain recurrence class for such a diffeomorphism, then, $H=M$.
\end{prop}

With our results we are able to deduce some stronger statements for bi-Lyapunov (and Lyapunov) stable homoclinic classes for generic diffeomorphisms in manifolds of low dimensions.

In \cite{PotS}, our results allowed us to deduce the conjecture in dimension 2 since a homoclinic class with nonempty interior must admit a dominated splitting (this can be deduced from \cite{BDP} in the case the homoclinic class has nonempty interior), so, combining Theorem \ref{descdom} with Theorem \ref{maintheorem} we get the following Theorem.

\begin{teo}\label{genericsurface}
Let $f$ be a $C^1$-generic surface diffeomorphism having a bi-Lyapunov stable homoclinic
class. Then $f$ is conjugated to a linear Anosov diffeomorphism in $\T^2$.
\end{teo}

For the reader's convenience, in section \ref{SectionBajas}, we give a proof of Theorem \ref{genericsurface} which is independent from Theorem \ref{descdom}. It contains the ideas of Theorem \ref{descdom} but less technicalities.

Theorem \ref{descdomLYAPUNOV}, together with results of \cite{PujSam} (see also \cite{ABCD}), has the following interesting consequence in surface dynamics (we shall give more precise definitions, together with some direct corollaries in section \ref{SectionBajas}):

\begin{teo}\label{disipativosuperficie} Let $f$ be a $C^1$-generic surface diffeomorphism having a Lyapunov stable homoclinic class $H$ which has a dissipative periodic point. Then, $H$ is a hyperbolic attractor.
\end{teo}

Theorem \ref{descdom} together with Theorem \ref{maintheorem} combine to give the following interesting result which supports the extension of the conjecture to bi-Lyapunov stable classes and solves completely Problem 5.1 in \cite{ABD} for 3 dimensional manifolds.

\begin{prop}\label{dimension3} Let $H$ be a bi-Lyapunov stable homoclinic class for a $C^1$-generic diffeomorphism in dimension 3. Then, $H$ has nonempty interior.
\end{prop}

\begin{obs} The proof of this Proposition also gives us that: if $H$ is a bi-Lyapunov stable homoclinic class for a $C^1$-generic diffeomorphism in any dimension, which admits a codimension one dominated splitting, $T_HM=E\oplus F$ with $\dim F=1$, and a periodic point $p$ of index $d-1$, then, the class has nonempty interior. \finobs
\end{obs}

\subsection{Idea of the proof}

The idea of the proof of Theorem \ref{maintheorem} is the following.

First we prove that if the homoclinic class is bi-Lyapunov stable, the
periodic points in the class (which are all saddles) should have
eigenvalues (in the $F$ direction) exponentially (with the period)
far from $1$. Otherwise we manage to obtain a sink or a source
inside the class (by using an improved version of Franks' lemma by Gourmelon, \cite{Go3}, which allows to
perturb the derivative controlling the invariant manifolds and also using Lyapunov stability) which is a contradiction.
After this is done, one can conclude with the same arguments as in \cite{PotS}.

To get the dominated splitting in the class (Theorem \ref{descdom}) we use a similar idea, if the class does not admit any dominated splitting, then, by \cite{BDP} one can create a sink or a source by perturbation, so, using again the Franks' lemma with control of the invariant manifolds, we are able to ensure that the sink or source is contained in the class using also Lyapunov stability and thus reach a contradiction. A similar idea also gives information on the index of the dominated splitting.

 To be able to use the improved version of Franks' lemma given by Gourmelon, we shall use a recent result by Bochi and Bonatti (\cite{BoBo}) which we present in section \ref{sectionsurface}.


\subsection{Organization of the paper}

 In section \ref{sectionmaintheorem} we shall prove Theorem \ref{maintheorem} and in section \ref{sectionsurface} we shall prove Theorem \ref{descdom}. The main argument used in the paper is presented with full details in the end of the proof of Lemma \ref{normadiferencial}, the other places where we use it will refer to it.
 In section \ref{SectionBajas} we prove the results in low dimensions, namely Theorem \ref{genericsurface} (with a proof independent from Theorem \ref{descdom}), Theorem \ref{disipativosuperficie} and Proposition \ref{dimension3} using Theorem \ref{maintheorem}.

  We present  the proof of Proposition \ref{lejosdetangencias} based on a Theorem of Yang \cite{Y} in section \ref{sectionY}. We present also a proof of Yang's theorem in Appendix B based in \cite{C}.

  In Appendix A we recall some generic results which are very well known and shall be used in the course of the paper.

\medskip

\textit{Acknowledgements:} Many of the ideas here came out from conversations with Martin Sambarino, I would like to thank him for that and for reading the first version and making suggestions. Discussions with Nikolaz Gourmelon were useful to improve the presentation and better understand his results which are the key tool used in this paper. I would not have been able to write this paper without the help of Sylvain Crovisier, I could explicitly detail his contributions to this paper, but it would be too long, thanks for all, but specially for the patience. The anonymous referee certainly gave an exhaustive read to the draft and made several important suggestions, specially finding some errors in some proofs and improving the presentation, I would like to thank him/her for that.

\section{Proof of Theorem \ref{maintheorem}}\label{sectionmaintheorem}

Let $Per(f)$ denote the set of periodic points of $f$. For $p\in Per(f)$, $\pi(p)$ denotes the period of $p$. We denote as $Per_{\alpha}(f)$ the set of index $\alpha$ periodic points. Let $\mathcal{O}$ be a periodic orbit and $E$ a $Df$ invariant subbundle of $T_{\mathcal O}M$; $D_{\mathcal O}f_{/E}$ denotes the cocycle over the periodic orbit given by its derivative restricted to the invariant subbundle.

We shall make some definitions. Let $\mathcal O$ be a periodic orbit and $\cA_\cO$ be a linear cocycle \footnote{A linear cocycle $\cA$ of dimension $n$ over a transformation $f:\Sigma \to \Sigma$ is a map $A:\Sigma \to GL(n,\R)$. When one point $p\in \Sigma$ is $f-$periodic, the eigenvalues of the cocycle at $p$ are the eigenvalues of the  matrix given by $A_{f^{\pi(p)-1}(p)} \ldots A_p$. See \cite{BGV}.} over $\mathcal O$. We say that $\cA_{\cO}$ has a \emph{strong stable manifold of dimension} $i$ if the eigenvalues $|\lambda_1|\leq |\lambda_2| \leq \ldots \leq |\lambda_d|$ of $\cA_\cO$ satisfy that $|\lambda_i| < \min\{1,|\lambda_{i+1}|\}$. If the derivative of $\mathcal O$ has strong stable manifold of dimension $i$ then classical results ensure the existence of a local, invariant manifold tangent to the the subspace generated by the eigenvectors of these $i$ eigenvalues and imitating the behavior of the derivative (see \cite{HPS}).

 Let $\Gamma_i$ be the set of cocycles over $\mathcal O$ which have a strong stable manifold of dimension $i$.

We endow $\Gamma_i$ with the following distance, $d(\cA_\cO,\cB_\cO)=\max\{ \|\cA_\cO - \cB_\cO\| , \|\cA_\cO^{-1} - \cB_\cO^{-1}\|\}$ where the norm is $\|\cA_\cO\|=\sup_{p\in \mathcal O } \{ \frac{\|A_p(v)\|}{\|v\|} \ ; \ v \in T_pM\backslash \{0\}\}$.

Let $g$ be a perturbation of $f$ such that the cocycles $D_{\mathcal O}f$ and $D_{\mathcal O} g$ are both in $\Gamma_i$, and let $U$ be a neighborhood of $\mathcal O$. We shall say that $g$ \emph{preserves locally the $i-$ strong stable manifold of $f$ outside $U$}, if the set of points of the $i-$strong stable manifold of $\mathcal O$ outside $U$ whose positive iterates do not leave $U$ once they entered it, are the same for $f$ and for $g$.

We have the following theorem due to Gourmelon which will be the key tool for proving  the results here presented.

\begin{teo}[\cite{Go3}]\label{frankslemmaGourmelon}
Let $f$ be a diffeomorphism, and $\mathcal{O}$ a periodic orbit of $f$ such that $D_{\mathcal O}f \in \Gamma_i$ and let $\gamma:[0,1]\to \Gamma_i$ be a path starting at $D_{\mathcal O}f$. Then, given a neighborhood $U$ of $\mathcal O$, there is a perturbation $g$ of $f$ such that $D_{\mathcal O}g=\gamma(1)$, $g$ coincides with $f$ outside $U$ and preserves locally the $i-$strong stable manifold of $f$ outside $U$. Moreover, given $\U$ a $C^1$ neighborhood of $f$, there exists $\eps>0$ such that if $\diametro (\gamma)<\eps$ one can choose $g \in \U$.
\end{teo}

We observe that Franks' lemma (see \cite{frankslema}) is the previous theorem with $i=0$. Also, we remark that Gourmelon's result is more general since it allows to preserve more than one stable and more than one unstable manifolds (of different dimensions, see \cite{Go3}).

\begin{lema}\label{normadiferencial} Let $H$ be a homoclinic class
which is Lyapunov stable of a $C^1$-generic diffeomorphism $f$ such that the class has only periodic orbits of index smaller or equal to $\alpha$. So, there exists $K_0>0$, $\lambda \in (0,1)$ and $m_0 \in \Z$ such that for every $p\in Per_{\alpha}(f_{/H})$ of sufficiently large period one has

$$\prod_{i=0}^k \left\|\prod_{j=0}^{m_0-1} Df^{-1}_{/E^u(f^{-im_0 -j}(p))}\right\| < K_0 \lambda^k \qquad k=\left[\frac{\pi(p)}{m_0}\right]. $$
\end{lema}

\dem{\!\!} Let $\Res$ be a residual subset of $\diff 1$ such that if $f\in \Res$ and $H$ is a Lyapunov stable homoclinic class of a periodic point $q$ of index $\alpha$, there exists a small neighborhood $\U$ of $f$ where the continuation $q_g$ of $q$ is well defined and such that for every $g\in \U \cap \Res$ one has that $H(g,q_q)$ is Lyapunov stable, and such that $g$ is a continuity point of the map $g\mapsto H(g,q_g)$. Also, being $f$ generic, we can assume that for every $g \in \U \cap \Res$ and every $p\in Per_{\alpha}(g)\cap H(g,q_g)$ we have that $H(g,q_g)=H(g,p)$, so, the orbits of $p$ and $q_g$ are homoclinically related (see the Appendix A, we are using properties $a1)$, $a2)$, $a3)$, $a4)$ and $a5)$ of Theorem A.1).

We can also assume that $\U$ and $\Res$ were chosen so that for every $g \in \U \cap \Res$ every periodic point in $H(g,q_g)$ has index smaller or equal to $\alpha$. We can also assume that $q_g$ has index $\alpha$ for every $g \in \U$ (see property $a5)$ of Theorem A.1).

Lemma II.5 of \cite{manie1} asserts that to prove the thesis is enough to show that there exists $\eps>0$ such that the set of cocycles $\Theta_{\alpha}=\{ D_{\mathcal O(p)}f^{-1}_{/E^u} \ :  \ ,
\ p\in Per_{\alpha}(f_{/H}) \}$ which all have its eigenvalues of modulus bigger than one, verify that every $\eps-$perturbation of them preserves this property. That is, given $p\in Per_\alpha (f_{/H})$ one has that every $\eps-$perturbation $\{A_0,\ldots, A_{\pi(p)-1}\}$ of $D_{\mathcal O(p)}f$ verifies that $A_{\pi(p)-1}\ldots A_0$ has all its eigenvalues of modulus bigger or equal to one.

Therefore, assuming by contradiction that the Lemma is false, we get that $\forall \eps>0$ there exists a periodic point $p\in Per_\alpha(f_{/H})$ and a linear cocycle over $p$, $\{A_0, \ldots, A_{\pi(p)}\}$ satisfying that  $\|D_{f^i(p)}f_{/E^u} -A_i\|\leq \eps$ and  $\|D_{f^i(p)}f^{-1}_{/E^u} -A_i^{-1}\|\leq \eps$ and such that $\prod_{i=0}^{\pi(p)-1}A_i$ has some eigenvalue of modulus smaller or equal to $1$.

In coordinates $T_{\mathcal O(p)}M =E^u\oplus (E^u)^\perp$, since $E^u$ is invariant we have that the form of $Df$ is given by

\[ D_{f^i(p)}f= \left(
     \begin{array}{cc}
       D_{f^i(p)}f_{/E^u} & K^1_i(f) \\
       0 & K^2_i(f) \\
     \end{array}
   \right)
 \]

Let $\gamma:[0,1]\to \Gamma_{\alpha}$ given in coordinates $T_{\mathcal O(p)}M =E^u\oplus (E^u)^\perp$ by

$$\gamma_i(t)= \left(
     \begin{array}{cc}
      (1-t)D_{f^i(p)}f_{/E^u} + t A_i & K^1_i(f) \\
       0 & K^2_i(f) \\
     \end{array}
   \right)$$

whose diameter is bounded by $\eps$ (see Lemma 4.1 of \cite{BDP}).

\medskip

Now(\footnote{The following argument will be used repeatedly along the paper.}), choose a point $x$ of intersection between $W^s(p,f)$ with $W^u(q,f)$ and choose a neighborhood $U$ of the orbit of $p$ such that:

\begin{itemize}
\item[(i)] It does not intersect the orbit of $q$.
\item[(ii)] It does not intersect the past orbit of $x$.
\item[(iii)] It verifies that once the orbit of $x$ enters $U$ it stays there for all its future iterates by $f$.
\end{itemize}

It is very easy to choose $U$ satisfying (i) since both the orbit of $p$ and the one from $q$ are finite. Since the past orbit of $x$ accumulates in $q$ is not difficult to choose $U$ satisfying (ii). To satisfy (iii) one has only to use the fact that $x$ belongs to the stable manifold of $p$ so, after a finite number of iterates, $x$ will stay in the local stable manifold of $p$, and it is not difficult then to choose a neighborhood $U$ which satisfies (iii) also.

 Applying Theorem \ref{frankslemmaGourmelon} we can perturb $f$ to a new diffeomorphism $\hat g$ so that the orbit of $p$ has index greater than $\alpha$ and so that it preserves locally its strong stable manifold. This allows to ensure the intersection between $W^u(q_{\hat{g}},\hat g)$ and $W^s(p,\hat{g})$.

This intersection is transversal so it persist by small perturbations, the same occurs with the index of $p$ so we can assume that $\hat g$ is in $\Res \cap \U$. Using Lyapunov stability of $H(\hat g, q_{\hat g})$ we obtain that $p \in H(\hat g, q_{\hat g})$. This is because a Lyapunov stable homoclinic class is saturated by unstable sets, so, since $q_{\hat g} \en H(\hat{g},q_{\hat{g}})$ we have that $\overline{W^u(\hat{g},q_{\hat{g}})} \en H(\hat{g},q_{\hat{g}})$ and since $W^s(p,\hat g)\cap W^u(\hat{g},q_{\hat{g}}) \neq \emptyset$, we get that $p \in \overline{W^u(\hat{g},q_{\hat{g}})}$ and we get what we claimed.

This contradicts the choice of $\mathcal U$  since we find a diffeomorphism in $\U\cap \Res$ with a periodic point with index bigger than $\alpha$ in the continuation of $H$, and so the lemma is proved.

\lqqd

\begin{obs}\label{hipeneelperiodocodimensionuno} One can recover Lemma 2 of \cite{PotS} in this context. In fact, if there is a codimension one dominated splitting of the form $T_HM =E\oplus F$ with $\dim F=1$ then (using the adapted metric given by \cite{GourmelonAdaptada}) for a periodic point of maximal index one has $\|Df^{-1}_{/F(p)}\|\leq \|Df^{-1}_{/E^u(p)}\|$ so,
$$\prod_{i=0}^k \left\|\prod_{j=0}^{m_0-1} Df^{-1}_{/F(f^{-im_0 -j}(p))}\right\| < K_0 \lambda^k \qquad k=\left[\frac{\pi(p)}{m_0}\right] $$

And since $F$ is one dimensional one has $\prod_i\|A_i\| =\|\prod_i A_i\|$ so $\|Df^{-\pi(p)}_{/F(p)}\|\leq K_0 \lambda^{\pi(p)}$ (maybe changing the constants $K_0$ and $\lambda$).

In fact, there is $\gamma \in (0,1)$ such that for every periodic point of maximal index and big enough period one has $ \|Df^{-\pi(p)}_{/F(p)}\|\leq \gamma^{\pi(p)}$

Also, it is not hard to see, that if the class admits a dominated splitting of index bigger or equal than the index of all the periodic points in the class, then, periodic points should be hyperbolic in the period along $F$ (for a precise definition and discussion on this topics one can read \cite{BGY}, \cite{wen}).
\finobs
\end{obs}

\begin{obs}\label{deapedazos} As a consequence of the proof of the lemma we get that: One can perturb the eigenvalues along an invariant subspace of a cocycle  without altering the rest of the eigenvalues. The perturbation will be of similar size to the size of the perturbation in the invariant subspace. See Lemma 4.1 of \cite{BDP}. Notice also that we could have perturbed the cocycle $\{K_i^2(f)\}_i$ without altering the eigenvalues of the cocycle $D_{\mathcal O(p)}f_{/E^u}$.\finobs
\end{obs}

One can now conclude the proof of Theorem \ref{maintheorem} with the same techniques as in the proof of the main Theorem of \cite{PotS}.

We have that $T_HM=E\oplus F$ with $\dim F=1$. We first prove that the center unstable curves tangent to $F$ should be unstable and with uniform size (this is Lemma 3 of \cite{PotS}). To do this, we first use Lemma \ref{normadiferencial} to get this property in the periodic points and then use the results from \cite{PS2} and \cite{BC} to show that the property extends to the rest of the points. This dynamical properties imply also uniqueness of these central unstable curves.

Assuming the bundle $F$ is not uniformly expanded, one has two cases, one can apply Liao's selecting lemma or not (see \cite{wen}).

In the first case one gets weak periodic points inside the class which contradict the thesis of Lemma \ref{normadiferencial}. The second case is similar, if Liao's selecting lemma does not apply, one gets a minimal set inside $H$ where $E$ is uniformly contracting and, using the dynamical properties of the center unstable curves, classical arguments give a periodic point close to this minimal set. Since the stable manifold of this periodic point will be uniform, it will intersect the unstable manifold of a point in $H$, and then Lyapunov stability implies the point is inside the class and again contradicts Lemma \ref{normadiferencial}.

For more details see \cite{PotS}.

\section{Low dimensional consequences}\label{SectionBajas}

We shall first prove Theorem \ref{genericsurface}. We remark that it is an easy corollary of Theorems \ref{descdom} and \ref{maintheorem} since together they give hyperbolicity of the homoclinic class. However, this proof does not use Theorem \ref{descdom} and we believe it may illuminate the idea of the more general proof.

\begin{teo*} Let $H$ be a bi-Lyapunov stable homoclinic class of a $C^1$-generic surface diffeomorphism $f$. Then, $H$ admits a dominated splitting and thus $H=\T^2$ and $f$ is an Anosov diffeomorphism.
\end{teo*}

\dem{\!\!} Consider a generic diffeomorphism $f$ satisfying properties $a1)$, $a2)$, $a3)$ and $a4)$ of Theorem A.1. Consider a periodic point $q\in H$ fixed such that for a neighborhood $\mathcal U$ of $f$ the class $H(q_g,g)$ is bi-Lyapunov stable for every $g$ in a residual subset of $\mathcal U$.

If the class does not admit a dominated splitting, arguments from \cite{manie1} allow us to deduce that after a small perturbation, for any $\delta>0$ we have a periodic point $p$ such that the angle between the stable and unstable bundle is smaller than $\delta$. Let us assume, for example, that the determinant of $Df_p^{\pi(p)}$ is smaller than one.

By classical arguments (see for example Lemma 7.7 of \cite{beyondhip}) one knows that after composing $D_{\mathcal O(p)}f$ with a rotation of angle smaller than $\delta$ one has that the resulting cocycle, $\cA_{\cO(p)}$, has complex eigenvalues. These eigenvalues have modulus smaller than $1$ (since the rotation has determinant equal to $1$). Also, the same argument of the determinant implies that we can join $D_{\mathcal O(p)}f$ with $\cA_{\cO(p)}$ by a curve with small diameter and such that every cocycle in the curve has one eigenvalue with modulus smaller than one. In fact, by stopping the curve before getting both eigenvalues equal, one is assured to be in the hypothesis of Theorem \ref{frankslemmaGourmelon}.

We can then apply Theorem \ref{frankslemmaGourmelon} to perturb the periodic point $p$ preserving its strong stable manifold locally.

Arguing as in Lemma \ref{normadiferencial} we can do the perturbation so that $W^{u}(q_g,g)$ intersects $W^s(p,g)$ and thus, using Lyapunov stability we obtain that $p \in H(q_g)$ which is absurd since $p$ is a sink.

The rest follows by applying Theorem \ref{maintheorem} which implies that $H$ is hyperbolic, and thus, using local product structure and bi-Lyapunov stability, we get that $H$ is open and closed, and thus $H=M$. So, $f$ is an Anosov diffeomorphism, and by Franks' theorem (\cite{FranksAnosov}), it must be conjugated to a linear Anosov diffeomorphism of $\T^2$.

\lqqd

We say that a periodic point $p$ of $f$ is \emph{dissipative} if the determinant of $Df_p^{\pi(p)}$ is smaller than one. We can now prove Theorem \ref{disipativosuperficie}

\begin{teo*} Let $f$ be a $C^1$-generic surface diffeomorphism having a Lyapunov stable homoclinic class $H$ which has a dissipative periodic point. Then, $H$ is a hyperbolic attractor.
\end{teo*}

\dem{\!\!} Let $f$ be a $C^1$-generic surface diffeomorphism satisfying properties $a1)$, $a2)$, $a3)$ and $a4)$ of Theorem A.1. We also assume that $f$ satisfies another generic property given in Theorem 2 of \cite{ABCD} (see also \cite{PujSam}) stating that a chain recurrence class admitting a dominated splitting is hyperbolic (this holds only in surfaces).

The argument above gives that the homoclinic class has a dominated splitting (notice that we assumed that there was a periodic point with determinant $<1$ and only used Lyapunov stability). The above considerations give that the homoclinic class is hyperbolic, and since it is Lyapunov stable it is a hyperbolic attractor.
\lqqd

\begin{obs} Notice that for Theorem \ref{genericsurface} we did not use the results of \cite{PujSam} which involve $C^2$ approximations. In the Theorem above, we have only used Lyapunov stability for $f$, and so Theorem \ref{maintheorem} allows only to prove hyperbolicity of one of the bundles, that is why we need to use this new generic property given by Theorem 2 of \cite{ABCD}.\finobs
\end{obs}

The last theorem has some immediate consequences which may have some interest on their own.

We say that an embedding $f:\D^2 \to \D^2$ is \emph{dissipative} if for every $x\in \D^2$ we have that $|det(D_xf)|<b<1$. Recall that for a dissipative embeddings of the disc, the only hyperbolic attractors are the sinks (\cite{Plykin}).

\begin{cor} Let $f:\D^2 \to \D^2$ be a generic dissipative embedding. Then, every Lyapunov stable homoclinic class is a sink.
\end{cor}

This gives in particular, that for the well known Henon map (see \cite{beyondhip} chapter 4), $C^1$-generic diffeomorphisms near, do not have a strange attractors (notice that they may have aperiodic quasi-attractors, which is not hard to show would be semiconjugated to adding machines, see \cite{beyondhip} chapter 10).

Now, we shall prove Proposition \ref{dimension3}. This gives that Question \ref{Question} has the same answer as Conjecture \ref{Conjetura} for diffeomorphisms in $3$ dimensional manifolds. Compare with the results from \cite{PotS}.

\dem{ of Proposition \ref{dimension3}} Applying Theorem \ref{descdom} one can assume that the class $H$ admits a dominated splitting of the form $E\oplus F$, and without loss of generality one can assume that $\dim F=1$.

Theorem \ref{maintheorem} thus implies that $F$ is uniformly expanded so the splitting is $T_HM = E \oplus E^u$.

Assume first that there exist a periodic point $p$ in $H$ of index $2$. Thus, this periodic point has a local stable manifold of dimension $2$ which is homeomorphic to a 2 dimensional disc.

Since the class is Lyapunov stable for $f^{-1}$ the stable manifold of the periodic point is completely contained in the class.

Now, using Lyapunov stability for $f$ and the foliation by strong unstable manifolds given by \cite{HPS} one gets (saturating by unstable sets the local stable manifold of $p$) that the homoclinic class contains an open set. This implies the thesis under this assumption.

So, we must show that if all the periodic points in the class have index $1$ then the class is the whole manifold. As we have been doing, using the genericity of $f$ we can assume that there is a residual subset $\Res$ of $\diff 1$ and an open set $\U$ of $f$ such that for every $g\in \U \cap \Res$ all the periodic points in the class have index $1$.

We have 2 situations, on the one hand, we consider the case where $E$ admits two invariant subbundles, $E=E^1\oplus E^2$, with a dominated splitting and thus, we get that $E^1$ should be uniformly contracting (using Theorem \ref{maintheorem}) proving that the homoclinic class is the whole manifold (Corollary \ref{Cor1}).

If $E$ admits no invariant subbundles then, using Theorem \ref{BGV} below, we can perturb the derivative of a periodic point in the class, so that the cocycle over the periodic point restricted to $E$ has all its eigenvalues contracting. So, we can construct a periodic point of index $2$ inside the class.\lqqd

\begin{obs} It is very easy to adapt the proof of this proposition to get that: If a bi-Lyapunov stable homoclinic class of a generic diffeomorphism admits a codimension one dominated splitting, $T_HM=E\oplus F$ with $\dim F=1$, and has a periodic point of index $d-1$, then, the class has nonempty interior.\finobs
\end{obs}

\section{Existence of a dominated splitting}\label{sectionsurface}

We  prove here Theorem \ref{descdom} and \ref{descdomLYAPUNOV} which state that a bi-Lyapunov stable homoclinic class of a generic diffeomorphism (or a Lyapunov stable homoclinic class with a dissipative periodic orbit) admits a dominated splitting.

The idea is the following: in case $H$ does not admit any dominated splitting we can perturb the derivative of some periodic point in order to convert it into a sink or a source with the techniques of \cite{BDP} and \cite{BGV}. We pretend to use Theorem \ref{frankslemmaGourmelon} to ensure that the stable or unstable manifold of a periodic point in the class intersects the unstable or stable set of the source or the sink respectively and reach a contradiction. 

This technique generalizes to the case where the class does not admit a dominated splitting with index between the indices of the periodic points in the class since by using these kind of results we may construct either a periodic point of smaller index than those on the class or one of bigger index and manage to relate it to the points in the class using Lyapunov stability.

The section is divided in two, in the first part, we recall some of the notions defined in \cite{BGV} and we state a recent result by Bochi and Bonatti (see\footnote{I would like to thank C. Bonatti for providing me a preliminary version of \cite{BoBo}.} \cite{BoBo}) which will allow us to use Theorem \ref{frankslemmaGourmelon}. We give also an idea of the proof since it is a key step in the proof.

In the second part, we prove Theorem \ref{descdom} and Theorem \ref{descdomLYAPUNOV}.

\subsection{Perturbations of cocycles over small paths}
The goal of this section is to prove Theorem \ref{BGV} which allows to perturb linear cocycles without dominated splitting in order to use Gourmelon's Theorem \ref{frankslemmaGourmelon}.

Before we proceed with the statement  of Theorem \ref{BGV} we shall give some definitions taken from \cite{BGV} and others which we shall adapt to fit our needs.

Let $\mathcal{A} = (\Sigma, f, E, A)$ be a \emph{large period linear cocycle}\footnote{Sometimes, we shall abuse notation and call it just \emph{cocycle}.} of dimension $d$ bounded by $K$ over an infinite set $\Sigma$, that is

\begin{itemize}

\item[-] $f:\Sigma \to \Sigma$ is a bijection such that all points in $\Sigma$ are periodic and such that given $n>0$ there are only finitely many with period less than $n$.
\item[-]  $E$ is a vector bundle over $\Sigma$, that is, there is $p: E \to \Sigma$ such that $E_x=p^{-1}(x)$ is a vector space of dimension $d$ endowed with an euclidian metric $\langle,\rangle_x$.
\item[-] $A:x\in \Sigma \mapsto A_x\in GL(E_x, E_{f(x)})$ is such that $\|A_x\|\leq K$ and $\|A_x^{-1}\|\leq K$.

\end{itemize}

In general, we shall denote $A^\ell_x = A_{f^{\ell-1}(x)}  \ldots  A_x$ where juxtaposition denotes the usual composition of linear transformations. 

For every $x\in \Sigma$ we denote by $\pi(x)$ its period and $M_x^A = A^{\pi(x)}_x$ which is a linear map in $GL(E_x,E_x)$ (which allows to study eigenvalues and eigenvectors).

For $1\leq j \leq d$

$$\sigma^{j}(x,\mathcal A) = \frac{log|\lambda_j|}{\pi(x)}$$

Where $\lambda_1,\ldots \lambda_d$ are the eigenvalues of $M_x^A$ in increasing order of modulus. As usual, we call $\sigma^j(x,\cA)$ the $j-$th Lyapunov exponent of $\cA$ at $x$.

Given an $f-$invariant subset $\Sigma'\en \Sigma$, we can always restrict the cocycle to the invariant set defining the cocycle $\cA|_{\Sigma'} = (f|_{\Sigma'}, \Sigma', E|_{\Sigma'}, A|_{\Sigma'})$.

We shall say that a subbundle $F\en E$ is \emph{invariant} if $\forall x \in \Sigma$ we have $A_x(F_x) =F_{f(x)}$. When there is an invariant subbundle, we can write the cocycle in coordinates $F\oplus F^{\perp}$ (notice that $F^\perp$ may not be invariant) as

\[    \left(
       \begin{array}{cc}
         \cA_F & C_F \\
         0 & \cA|F \\
       \end{array}
     \right)\]

Where $C_F$ is uniformly bounded. This induces two new cocycles: $\cA_F =(\Sigma, f, F, A|_{F})$ on $F$ (where $A|_F$ is the restriction of $A$ to $F$) and $\cA|F = (\Sigma, f, E|F, A|F)$ on $E|F \simeq F^\perp$ where $(A|F)_x \in GL((F_x)^\perp, (F_{f(x)})^\perp)$ is given by $p^2_{f(x)} \circ A_x$ where $p^2_x$ is the projection map from $E$ to $F^\perp$. Notice that changing only $\cA_F$ affects only the eigenvalues associated to $F$ and changing only $\cA|F$ affects only the rest of the eigenvalues, recall Remark \ref{deapedazos}. See section 4.1 of \cite{BDP} for more discussions on this decomposition.

If $\cA$ has two invariant subbundles $F$ and $G$, we shall say that $F$ is $\ell$-\emph{dominated} by $G$ (and denote it as $F\prec_{\ell} G$) on an invariant subset $\Sigma'\en \Sigma$ if for every $x\in \Sigma'$ and for every vectors $v\in F_x \backslash \{0\}$, $w\in G_x \backslash \{0\}$ one has

\[ \frac{\|A^\ell_x(v)\|}{\|v\|} \leq \frac{1}{2} \frac{\|A^\ell_x(w)\|}{\|w\|} \ . \]

We shall denote $F\prec G$ when there exists $\ell>0$ such that $F\prec_{\ell} G$.

If there exists complementary invariant subbundles $E=F\oplus G$  such that $F \prec G$ on a subset $\Sigma'\en \Sigma$, we shall say that $\cA$ admits a \emph{dominated splitting} on $\Sigma'$. 

As in \cite{BGV}, we shall say that $\cA$ is \emph{strictly without domination} if it is satisfied that always that $\cA$ admits a dominated splitting in a set $\Sigma'$ it is satisfied that $\Sigma'$ is finite.

Let $\Gamma_x$ be the set of cocycles over the orbit of $x$ with the distance (see also section \ref{sectionmaintheorem})

$$  d(\cA_x, \cB_x) = {\displaystyle \sup_{0\leq i < \pi(x), v \in E\backslash \{0\}}} \left\{ \frac{\|(A_{f^i(x)}-B_{f^i(x)})v\|}{\|v\|} , \frac{\|(A_{f^i(x)}^{-1}-B_{f^i(x)}^{-1})v\|}{\|v\|} \right\} $$

 and let $\Gamma_{\Sigma}$ (or $\Gamma_{\Sigma,0}$) the set of bounded large period linear cocycles over $\Sigma$. Given $\cA \in \Gamma_\Sigma$ we denote as $\cA_x \in \Gamma_x$ to the cocycle $\{A_x, \ldots, A_{f^{\pi(x)-1}(x)} \}$.

We recall from section \ref{sectionmaintheorem} that we say that the cocycle $\cA_x$ has \emph{strong stable manifold of dimension} $i$ if  $\sigma^i(x,\cA_x) < \min \{0,\sigma^{i+1}(x, \cA_x) \}$.

For $0\leq i \leq d$, let

$$\Gamma_{\Sigma,i} = \{ \cA \in \Gamma_{\Sigma} \ : \ \forall x \in \Sigma \ ; \ \cA_x \ has \ strong \ stable \ manifold \ of \ dimension \ i \}$$

\smallskip

Following \cite{BGV} we say that $\mathcal B$ is a \emph{perturbation} of $\mathcal A$ (denoted by $\mathcal B \sim \mathcal A$) if for every $\eps>0$ the set of points $x\in \Sigma$ such that $\cB_x$ is not $\eps-$close to the cocycle $\cA_x$ is finite.

Similarly, we say that $\mathcal B$ is a \emph{path perturbation} of $\mathcal A$ if for every $\eps>0$ one has that the set of points $x \in \Sigma$ such that $\cB_x$ is not a perturbation of $\mathcal A_x$ along a path of diameter $\leq \eps$ is finite. That is, there is a path $\tilde \gamma: [0,1] \to \Gamma_{\Sigma}$ such that $\tilde \gamma(0)=\cA$ and $\tilde \gamma(1) = \cB$ such that $\tilde \gamma_x:[0,1]\to \Gamma_x$ are continuous paths and given $\eps>0$ the set of $x$ such that $\tilde \gamma_x([0,1])$ has diameter $\geq \eps$ is finite.

In general, we shall be concerned with path perturbations which preserve the dimension of the strong stable manifold, so, we shall say that $\mathcal B$ is a \emph{path perturbation of index $i$} of a cocycle $\mathcal A \in \Gamma_{\Sigma,i}$ iff: $\cB$ is a path perturbation of $\cA$ and the whole path is contained in $\Gamma_{\Sigma,i}$.  This induces a relation in $\Gamma_{\Sigma,i}$ which we shall denote as $\sim^{\ast}_i$.

We have that $\sim$ and $\sim^\ast_i$ are equivalence relations in $\Gamma_\Sigma$ and $\Gamma_{\Sigma,i}$ respectively, and clearly $\sim^\ast_i$ is contained in $\sim$.

The \emph{Lyapunov diameter} of the cocycle $\mathcal A$ is defined as

$$\delta(\mathcal A)={\displaystyle \liminf_{\pi(x)\to \infty}} \  [\sigma^d(x,\mathcal A) - \sigma^1(x,\mathcal A)] .$$

If $\mathcal A \in \Gamma_{\Sigma,i}$, we define $\delta_{min}(\mathcal A) = \inf_{\mathcal B \sim \mathcal A} \{\delta(\mathcal B)\}$. Similarly, we define $\delta_{min}^{\ast,i}(\mathcal A)= \inf_{\mathcal B \sim^\ast_i \mathcal A} \{\delta(\mathcal B)\}$. Notice that $\delta_{min}^{\ast,i} (\mathcal A) \geq \delta_{min}(\cA)$ and a priori it could be strictly bigger.

\begin{obs}\label{dimensionextremal} For any cocycle $\cA$, it is easy to see that $\delta_{min}^{\ast,0}(\cA) = \delta_{min}(\cA)$. It sufficies to consider the path $(1-t)\cA + t\cB$ where $\cB$ is a perturbation of $\cA$ having the same determinant over any periodic orbit and verifying $\delta(\cB)=\delta_{min}(\cA)$ (see Lemma 4.3 of \cite{BGV} where it is shown that such a $\cB$ exists).
 \finobs
\end{obs}

We are now ready to state the following Theorem which will allow us to use the results from \cite{BGV} together with the improved version of Franks' Lemma given by Gourmelon. This Theorem follows directly from the results of \cite{BoBo} where a much more general result is proven. For completeness we shall make some comments on the proof which we shall omit in full generality.

\begin{teo}[Consequence of \cite{BoBo}]\label{BGV}
Let $\mathcal{A}=(\Sigma, f, E, A)$ be a bounded large period linear cocycle of dimension $d$. Assume that
\begin{itemize}
\item[-] $\cA$ is strictly without domination.
\item[-] $\mathcal{A} \in \Gamma_{\Sigma,i}$
\item[-] For every $x\in \Sigma$, we have $|det(M_x^A)|< 1$ (that is, for all $x\in \Sigma$ we have $\sum_{j=1}^d \sigma^j(x,\cA) <0$).
\end{itemize}
Then, $\delta_{min}^{\ast,i}(\cA)=0$. In particular, given $\eps>0$ there exists a point $x\in \Sigma$ and a path $\gamma_x$ of diameter smaller than $\eps$ such that $\gamma_x(0)=\cA_x$, the matrix $M_x^{\gamma_x(1)}$ has all its eigenvalues of modulus smaller than one, and such that $\gamma(t) \in \Gamma_{\Sigma,i}$ for every $t \in [0,1]$.
\end{teo}

It is well known since Ma\~ne (see \cite{manie1}) that this result holds in dimension 2, namely, a two dimensional cocycle strictly without domination admits a path perturbation which decreases the Lyapunov diameter even without affecting at all the determinant, in particular, if the determinant starts being smaller than one, there is a path perturbation of index $i$ which takes the Lyapunov diameter to zero (for any $i=0,1,2$).

In higher dimensions, a natural way to attack the problem is to use induction. A simple perturbation imitating Proposition 3.7 of \cite{BGV} allows to obtain a diagonal cocycle (that is, for every $x\in \Sigma$ all the Lyapunov exponents are different) which in turn give a lot of one dimensional invariant subbundles associated to each Lyapunov exponent, let us call them $E_1(x,\cA), \ldots, E_d(x, \cA)$.

This gives

\begin{lema}\label{valorespropiosdistintos} For every $\mathcal A \in \Gamma_{\Sigma,i}$, there exists $\mathcal B \sim^{\ast}_i \cA$ such that for every $x\in \Sigma$ the eigenvalues of $M_x^B$ have all different modulus and their modulus is arbitrarily near the original one in $M_x^A$, that is, $|\sigma^i(x,\cA)-\sigma^i(x,\cB)| \to 0$ as $\pi(x)\to \infty$ (in particular, $\delta(\cB)=\delta(\cA)$).
\end{lema}

When we wish to apply induction a problem appears, in general, an invariant subbundle of a strictly without domination cocycle needs not be strictly without domination (which would allow us to use induction). However, with more care, one can make a more subtle perturbation (which will be necessarily global) in order to preserve the existence of this invariant subbundle as well as breaking the domination and then be able to use induction.

Let us sketch the proof of the Theorem in the case of \cite{BGV} (without need of path perturbations) for three dimensional cocycles. First, we can assume that the cocycle is diagonal, that is, we have an invariant splitting of the form $E_1 \oplus E_2 \oplus E_3$ associated to the eigenvalues in increasing order of modulus. It is easy to show that we can consider the cocycle satisfying $\delta(\cA)=\delta_{min}(\cA)$. If $E_1$ is not dominated by $E_2$, with a small perturbation of $\cA_{E_1\oplus E_2}$ we can (using induction) reduce the Lyapunov diameter which contradicts $\delta(\cA)=\delta_{min}(\cA)$. So, we get that $E_1$ should be dominated by $E_2$. Since $\cA$ is strictly without domination, we get that $E_1|E_2$ cannot be dominated by $E_3|E_2$, so, again by a two dimensional perturbation we can decrease the Lyapunov diameter and conclude.

Notice that the last perturbation cannot be made if we wish to perform a small path perturbation of index $2$, since the eigenvalues associated to $E_3$ and $E_2$ may cross in this process loosing the dimension of the strong stable manifold.

To overcome this difficulty, we must make first a perturbation of $\cA_{E_2 \oplus E_3}$ in order to break the domination between $E_1$ and $E_2$. This is achieved by the following Proposition which we state without proof.

\begin{prop}\label{dim2tricotomia} Given $K>0$, $k>0$ and $\eps>0$, there exists $N>0$ and $\ell$ such that if
\begin{itemize}
\item[-] There exists a diagonal cocycle $\cA_x$ of dimension $2$ and bounded by $K$ over a periodic orbit of period $\pi(x)>N$.
\item[-] There exists a unit vector $v\in E_x$  such that

$$  \frac 3 2 \| A^\ell_{f^k(x)}|_{E_1}\| \geq \frac{ \|A^\ell_{f^k(x)} A^k_x v\|}{\|A^k_x v\|} $$

\end{itemize}

 Then, there exists $\cB_x$ a path perturbation of $\cA_x$ of diameter smaller than $\eps$ verifying that all along the path the cocycle has the same Lyapunov exponents and

$$    \| B^{k}_{x} |_{E_1(x,\cB)} \| \leq 2 \|A^{k}_{x}v\|$$

\end{prop}

Notice that the statement is quite general, in fact, some weaker statement could suffice. Roughly, the Proposition states that if the direction $E_1$ associated to the smaller eigenvalue has some iterates where there is a greater contraction in some vector $v$, then, one can perturb the cocycle along the orbit in a path, without changing the eigenvalues and obtaining a contraction similar to that of $v$ in the new eigendirection for some iterates of the cocycle. Notice that the important statement deals with the case where $v \notin E_1$ since if $v\in E_1$ there is no need for making a perturbation.  This allows to break domination as we shall explain later.

First, we give a quick idea of the proof of this proposition: The hypothesis of the Proposition gives us that we can send by a small perturbation along $\ell$ transformations that we can send the direction $v$ into the direction $E_1$. Since $E_2$ is more expanded than $E_1$ in the period, there will be some place where we can also send the direction of $E_1$ into the direction of $v$. This allows to show that the perturbation makes the $E_1$ direction to pass some iterates near the direction of $v$ and thus, inheriting its contraction.

To guaranty that the perturbation can be made in order two keep the eigenvalues one has to control the moment where the perturbations is made. Roughly, one has to show that one can take the direction $E_1$ into the $v$ direction sufficiently near the moment where the contraction takes place and that leaves time to ``correct'' the Lyapunov exponents along the rest of the orbit.

Now, we shall explain how to break the domination by using this Proposition. Consider a cocycle $\cA$ such that $F_j(\cA)= E_1(\cA) \oplus \ldots \oplus E_j(\cA)$ verifies that  $F_j \prec E_{j+1}(\cA)$ but $F_j(\cA)$ is not dominated by $E_{j+1}(\cA) \oplus E_{j+2}(\cA)$. We shall indicate how one can use the previous Proposition to find a new cocycle $\cB(\cA)$ which is a path perturbation of $\cA$ along a path that preserves all the Lyapunov exponents and verifies that $F_j(\cB)$ is not dominated by $E_{j+1}(\cB)$.

To perform this perturbation, we use the fact that since $F_j(\cA)$ is not dominated by $E_{j+1}(\cA) \oplus E_{j+2}(\cA)$ one can find points $x_n$ with periods converging to infinity such that there are vectors $v_n$ in $(E_{j+1}\oplus E_{j+2})_{x_n}$ which are more contracted along large periods of time than the vectors in $(F_{j})_{x_n}$. This allows (by breaking this time of strong contraction in two) us to use Proposition \ref{dim2tricotomia} to obtain a path perturbation $\cB$ which realizes this strong contraction along the bundle $(E_{j+1}(\cB))_{x_n}$ and thus breaking the domination as wanted.

After this is done, an induction argument permits to prove that if the subbundle $E_1(\cA) \oplus \ldots \oplus E_j(\cA)$ is strictly not dominated by the subbundle $E_{j+1}(\cA) \oplus \ldots \oplus E_{d}(\cA)$ then, there exists a path perturbation $\cB$ of $\cA$ such that the exponents are constant along the path and such that for $\cB$ the bundle $E_{j}(\cB)$ is strictly non dominated by $E_{j+1}(\cB)$.

This argument allows to translate the problem into a problem for two dimensional cocycles which as we said is already known how to deal with. With these observations made, one can conclude the proof of Theorem \ref{BGV}.

\subsection{Proof of Theorem \ref{descdom} and \ref{descdomLYAPUNOV}}

 Let $H$ be a generic bi-Lyapunov stable homoclinic class. Let us assume that $H$ contains periodic points of index $\alpha$ and we consider $\Delta^\eta_{\alpha}\en \{p \in Per_{\alpha}(f_{/H})\}$ the set of index $\alpha$ and $\eta-$disippative periodic points in $H$ for some $\eta<1$.

 It is enough to have one periodic point with determinant smaller than one to get that for some $\eta<1$, the set $\Delta^\eta_{\alpha}$ will be dense in $H$ (see Lemma 1.10 of \cite{BDP}). If no periodic point has determinant smaller than one, we work with $f^{-1}$ and use Lyapunov stability for $f^{-1}$ (see property $a1)$ in the Appendix A).

 Notice that if $H$ admits no dominated splitting, then neither does the cocycle of the derivatives over $\Delta^\eta_\alpha$. Then, we can assume that it is strictly without domination (maybe by considering an infinite subset). This implies that we can apply Theorem \ref{BGV} and there is a periodic point $p \in \Delta^\eta_{\alpha}$ which can be turned into a sink with a $C^1$ small perturbation done along a path contained in $\Gamma_{\Sigma,\alpha}$ (which maintains or increases the index).

Now we are able to use Theorem \ref{frankslemmaGourmelon} and reach a contradiction. Consider a periodic point $q \in \Delta^\eta_{\alpha}$ fixed such that for a neighborhood $\mathcal U$ of $f$ the class $H(q_g,g)$ is Lyapunov stable for every $g \in \U \cap \Res$ (recall the first paragraph of the proof of Lemma \ref{normadiferencial} and the Appendix A).

Suppose the class does not admit any dominated splitting, so, we have a periodic point $p\in \Delta^\eta_\alpha$ such that $f$ can be perturbed in an arbitrarily small neighborhood of $p$ to a sink for a diffeomorphism $g\in \U$ (which we can assume is in $\Res \cap \U$ since sinks are persistent) and preserving locally the strong stable manifold of $p$. So, we choose a neighborhood of $p$ such that it does not meet the orbit of $q$ nor the past orbit of some intersection of its unstable manifold with the local stable manifold of $p$ with the same argument as in Lemma \ref{normadiferencial}.

Thus, we get that $W^u(q_g,g) \cap W^s(p,g) \neq \emptyset$ and using Lyapunov stability we reach a contradiction since it implies that $p \in H(q_g)$ which is absurd since $p$ is a sink.

The same argument would work if $\Delta^\eta_{\alpha}$ consisted of points which where $\eta-$dissipative for  $f^{-1}$, but we should have used Lyapunov stability for $f^{-1}$ in that case.

Now we consider the finest dominated splitting of the class (that is, such that the subbundles admit no sub-dominated splitting, see \cite{beyondhip} Appendix B) which we denote as $T_HM = E_1 \oplus \ldots \oplus E_k$.

Assume that the indexes of the homoclinic class form the segment $[\alpha,\beta]$ (see property $a5)$ of Theorem A.1) and that there exists $l$ such that $\sum_{j=1}^{l-1} \dim E_j < \alpha$ and that $\sum_{j=1}^{l} \dim E_j >\beta$.

Then, Theorem \ref{BGV} allows to conclude that there exists a periodic point of index $\alpha$ such that an arbitrarily small perturbation of the cocycle $Df_{/E_l}$ has all its eigenvalues of the same modulus (notice that $E_l$ varies continuously, so, the cocycle $Df|_{E_l}$ is continuous). Assume that for example the modulus are smaller than one, so the previous argument gives also a contradiction by finding a periodic point in a perturbation of the class of index bigger or equal to $\sum_{j=1}^{l} \dim E_j$.

\lqqd

\begin{obs}
\begin{itemize}
\item Also the same ideas give that periodic points in the class must be volume hyperbolic in the period (not necessarily uniformly, see \cite{BGY} for a discussion on the difference between hyperbolicity in the period and uniform hyperbolicity)   for the extremal subbundles of the finest dominated splitting (see \cite{beyondhip} Appendix B for definitions).
\item Notice that the same proof works for Theorem \ref{descdomLYAPUNOV}. If a homoclinic class $H$ of a generic diffeomorphism $f$ is Lyapunov stable and the class has one periodic point $p$ such that $|det(Df^{\pi(p)}_p)| \leq 1$ then the class admits a dominated splitting.
    Moreover, this dominated splitting can be considered having index bigger or equal than the minimal index of the the class.

\item In fact, we can assume that if a Lyapunov stable homoclinic class admits no dominated splitting, then, there exists $\eta>1$ such that every periodic point $p$, it has determinant bigger than $\eta^{\pi(p)}$. Otherwise, there would exists a subsequence $p_n$ of periodic points with normalized determinant converging to $1$. After composing with a small homothety, we are in the hypothesis of Theorem \ref{descdomLYAPUNOV}.
\end{itemize}
\finobs
\end{obs}

\section{Bi-Lyapunov stable homoclinic classes far from tangencies}\label{sectionY}

We prove here that if $f$ is a generic diffeomorphism far from homoclinic tangencies and admits a chain recurrence class which is bi-Lyapunov stable, then $f$ must be transitive. First of all, one can reduce the study to homoclinic classes since in \cite{Y} it is proved that Lyapunov stable chain recurrence classes are in fact homoclinic classes.

For this we shall use Theorem \ref{maintheorem} and a recent result of \cite{Y}. For completeness and since we believe it may contribute to a better understanding, we present a proof of the theorem of \cite{Y} using new techniques introduced in \cite{C} in Appendix B of this paper.

First we state the following important result due to Yang for generic Lyapunov stable homoclinic classes far from tangencies (we shall denote $\overline{\Tang}$ as the set of diffeomorphisms which can be $C^1$ perturbed to get a homoclinic tangency between the stable and unstable manifold of a hyperbolic periodic point):

\begin{teo}[\cite{Y} Theorem 3]\label{teoY} Let $f \in \Res$, where $\Res$ is a residual subset of  $\diff 1 \backslash \overline{\Tang}$, and let $H$ be a Lyapunov stable homoclinic class for $f$ of minimal index $\alpha$. Let $T_HM=E\oplus F$ be a dominated splitting for $H$ with $\dim E=\alpha$, so, one of the following two options holds:

\begin{enumerate}
\item $E$ is uniformly contracting.
\item $E$ decomposes as $E^s\oplus E^c$ where $E^s$ is uniformly contracting and $E^c$ is one dimensional and $H$ is the Hausdorff limit of periodic orbits of index $\alpha-1$.
\end{enumerate}
\end{teo}

With this theorem and Theorem \ref{maintheorem}, Proposition \ref{lejosdetangencias} follows easily. We must remark that this theorem alone is enough to prove the same result for homoclinic classes with interior since they are not compatible with being accumulated in the Hausdorff topology with periodic points of index smaller than the class. Theorem \ref{maintheorem} is necessary in the bi-Lyapunov stable case as will be seen clearly in the proof.

\dem{of Proposition \ref{lejosdetangencias}} First of all, if the class has all its periodic points with index between $\alpha$ and $\beta$ we know that it admits a 3 ways dominated splitting of the form $T_HM=E\oplus G \oplus F$ where $\dim E=\alpha$ and $\dim F= d-\beta$. This is because we can apply the result of \cite{W} which says that far from homoclinic tangencies there is an index $i$ dominated splitting over the closure of the index $i$ periodic points together with the fact that index $\alpha$ and $\beta$ periodic points should be dense in the class since the diffeomorphism is generic (see Theorem A.1).

Now, we will show that $H$ admits a strong partially hyperbolic splitting. If $E$ is one dimensional, then it must be uniformly hyperbolic because of Theorem \ref{maintheorem}. If not, suppose $\dim E >1$ then, if it is not uniform, Theorem \ref{teoY} implies that it can be decomposed as a uniform bundle together with a one dimensional central bundle, since $\dim E>1$ we get a uniform bundle of positive dimension.

The same argument applies for $F$ using Lyapunov stability for $f^{-1}$ so we get a strong partially hyperbolic splitting.

Corollary 1 of \cite{ABD} finishes the proof.

\lqqd

\section*{Appendix A. Some generic properties}\label{sectionappendix}

The following theorem gives some well known generic properties we shall be using in the proof of our results. The main reference is \cite{beyondhip} (and the references therein) and mainly \cite{CMP}, \cite{C3}, \cite{BC} and \cite{ABCDW}.

\begin{teoA1}
 There exists a residual subset $\Res$ of $\diff 1$ such that if $f \in \Res$
\begin{itemize}
 \item[a1)] $f$ is Kupka-Smale (that is, all its periodic
  points are hyperbolic and their invariant manifolds intersect transversally). Also, we can assume that the determinant of the derivative of every periodic point is different from $1$.
 \item[a2)] The periodic points of $f$ are dense in the chain
 recurrent set of $f$. Moreover, if a chain recurrence class $C$ contains a periodic point $p$ then $C=H(p)$.

  \item[a3)] Given $p\in Per(f)$ there exists $\U_1$ a neighborhood of $f$ such that for every $g\in \U_1\cap \Res$, $g$ is a continuity point for the map $g \mapsto H(g,p_g)=H_g$ where $p_g$ is the continuation of $p$ for $g$. The continuity is with respect to the Hausdorff distance between compact subsets of $M$.

 \item[a4)] If a homoclinic class $H$ is Lyapunov stable for $f\in \Res$ then there exists $\U_2 \en \U_1$ a neighborhood of $f$ such that $H_g$ is Lyapunov stable for every $g\in \U_2\cap \Res$.

 \item[a5)] If $f\in \Res$ and $H$ is a homoclinic class. There exists an interval $[\alpha,\beta]$ of natural numbers such that for every $g\in \U_3 \en \U_1$, $H_g$ has periodic points of every index in $[\alpha,\beta]$ and every periodic point in $H_g$ has its index in that interval. Also, all periodic points of the same index are homoclinically related.

  \item[a6)] Let $X\en M$, and $x,y\in X$ such that for every $\eps>0$ there exists an $\eps-$pseudo orbit inside $X$ from $x$ to $y$. Then, for every $\delta>0$ there exists a segment of orbit $z,\ldots f^n(z)$ such that is contained in an $\delta$ neighborhood of $X$ and such that $d(x,z)<\delta$ and $d(y,f^n(z))<\delta$. Moreover, if $X$ is chain transitive, it is the Haussdorff limit of periodic orbits.
\end{itemize}
\end{teoA1}

We also recall the following perturbation result which will be used in the next appendix, it is called the Hayashi's connecting lemma, the statement we give is Theorem 5 of \cite{C3}:

\begin{teoA2}
Let $f\in \Diff^1(M)$ and $\U$ a neighborhood of $f$. Then, there exists $N>0$ such that every non periodic point $x \in M$ which admits two neighborhoods $W\en \hat W$ with the following property: for every $p,q \in M \backslash (f(\hat W) \cup \ldots \cup f^{N-1}(\hat W))$ such that $p$ has a forward iterate $f^{n_p}(p)\in W$ and $q$ has a backward iterate $f^{-n_q}(q)\in W$, we have that there exists $g \in \mathcal U$ which coincides with $f$ in $M \backslash (f(\hat W) \cup \ldots \cup f^{N-1}(\hat W))$ and such that for some $m>0$ we have $g^m(p)=q$.

Moreover, $\{p,g(p), \ldots, g^m(p)\}$ is contained in the union of the orbits $\{p, \ldots, f^{n_p}(p)\}$, $\{ f^{-n_q}(q), \ldots, q \}$ and the neighborhoods $\hat W, \ldots, f^N(\hat W)$. Also, the neighborhoods $\hat W, W$ can be chosen arbitrarily small.
\end{teoA2}

\section*{Appendix B. Lyapunov stable classes far from tangencies}\label{sectionYANG}

In this Appendix, we present a proof of Theorem \ref{teoY} (originally proved in \cite{Y}). We believe that having another approach to this important result is not entirely devoid of interest.

As in \cite{Y}, the proof has 3 stages, the first one is to reduce the problem to the central models introduced by Crovisier, the second one to treat the possible cases and finally, the introduction of some new generic property allowing to conclude in the difficult case.

Our proof resembles that of \cite{Y} in the middle stage (which is the most direct one after the deep results of Crovisier) and has small differences mainly in the other two.

For the first one, we use a recent result of \cite{C} (Theorem B.1) and for the last one, we introduce Lemma B.1 which can be compared with the main Lemma of \cite{Y} but the proof and the statement are somewhat different (in particular, ours is slightly stronger). We believe that this Lemma can find some applications (see for example \cite{C4}).

Before we start with the proof, we shall state the Theorem we shall use

\begin{teoB1}[Theorem 1 of \cite{C}]\label{splitting} Let $f \in \Res$ where $\Res$ is a residual subset of $\diff 1 \backslash \overline{\Tang}$ and $K_0$ an invariant compact set with dominated splitting $T_{K_0}M= E \oplus F$. If $E$ is not uniformly contracted, then, one of the following cases occurs.

\begin{enumerate}
\item $K_0$ intersects a homoclinic class whose minimal index is strictly less than $\dim E$.
\item $K_0$ intersects a homoclinic class whose minimal index is $\dim E$ and wich contains weak periodic orbits (for every $\delta$ there is a sequence of hyperbolic periodic orbits homoclinically related which converge in the Hausdorff topology to a set $K\en K_0$, whose index is $\dim E$ but whose maximal exponent in $E$ is in $(-\delta,0)$). Also, this implies that every homoclinic class $H$ intersecting $K_0$ verifies that it admits a dominated splitting of the form $T_HM = E'\oplus E^c \oplus F$ with $\dim E^c =1$.
\item There exists a compact invariant set $K\en K_0$ with minimal dynamics and which has a partially hyperbolic structure of the form $T_KM=E^s\oplus E^c\oplus E^u$ where $\dim E^c=1$ and $\dim E^s<\dim E$. Also, any measure supported on $K$ has zero Lyapunov exponent along $E^c$.
\end{enumerate}
\end{teoB1}

\dem{of Theorem \ref{teoY}} Let $\Res \en \diff 1 \backslash \overline{\Tang}$ be a residual subset such that for every $f \in \Res$ and every periodic point $p$ of $f$, there exists a neighborhood $\U$ of $f$, where the continuation $p_g$ of $p$ is well defined, such that $f$ is a continuity point of the map $g\mapsto H(g,p_g)$ and such that if $H(p,f)$ is a Lyapunov stable homoclinic class for $f$, then $H_g = H(p_g,g)$ is also Lyapunov stable for every $g\in \U\cap \Res$. Also, we can assume that for every $g\in \U\cap \Res$, the minimal index of $H_g$ is $\alpha$ (see Theorem A.1).

The class admits a dominated splitting of the form $T_HM=E \oplus F$ with $\dim E=\alpha$ (see \cite{W}). We assume that the subbundle $E$ is not uniformly contracted. This allows us to use Theorem B.1.

Since the minimal index of $H$ is $\alpha$, option $1)$ of the Theorem cannot occur.

We shall prove that option $3)$ implies option $2)$. That is, we shall prove that if $E$ is not uniformly contracted, then we are in option $2)$ of Theorem B.1.

This is enough to prove the Theorem since if we apply Theorem B.1 to $E'$ given by option $2)$ we get that since $\dim E' = \alpha-1$ option $1)$ and $2)$ cannot happen, and since option $3)$ implies option $2)$ we get that $E'$ must be uniformly contracted thus proving Theorem \ref{teoY} (observe that the statement on the Hausdorff convergence of periodic orbits to the class can be deduced from option $2)$ also by using Frank's Lemma).

\begin{claim} To get option $2)$ in Theorem B.1 is enough to find one periodic orbit of index $\alpha$ in $H$ which is weak (that is, it has one Lyapunov exponent in $(-\delta,0)$).
\end{claim}

\dem{\!\!}  This follows using the fact that being far from tangencies there is a dominated splitting in the orbit given by \cite{W2} with a one dimensional central bundle associated with the weak eigenvalue. Let $\mathcal O$ be the weak periodic orbit, so we have a dominated splitting of the form $T_{\mathcal O}M =E^s \oplus E^c \oplus E^u$.

  Using transitions (see \cite{BDP}), we can find a dense subset in the class of periodic orbits that spend most of the time near the orbit we found, say, for a small neighborhood $U$ of $\mathcal O$, we find a dense subset of periodic points $p_n$ such that the cardinal of the set $\{ i \in \Z \cap [0, \pi(p_n)-1] \ : \ f^i(p_n) \in U \}$ is bigger than $(1-\eps)\pi(p_n)$.

  Since we can choose $U$ to be arbitrarily small, we can choose $\eps$ so that the orbits of all $p_n$ admit the same dominated splitting (this can be done using cones for example\footnote{This argument is done carefully in the proof of Claim \ref{AfirmacionConos} in Subsection B.1.}) and maybe by taking $\eps$ smaller to show that $p_n$ are also weak periodic orbits.

\lqqd

  It rests to prove that option $3)$ implies the existence of weak periodic orbits in the class. To do this, we shall discuss depending on the structure of the partially hyperbolic splitting using the classification given in \cite{C}. There are 3 different cases. We shall name them as in the mentioned paper, but just explain the facts we will use.

  We have a compact invariant set $K\en H$ with minimal dynamics and which has a partially hyperbolic structure of the form $T_KM=E^s\oplus E^c\oplus E^u$ where $\dim E^c=1$ and $\dim E^s<\dim E$. Also, any measure supported on $K$ has zero Lyapunov exponent along $E^c$. We shall assume that the dimension of $E^s$ is minimal in the sense that every other compact invariant $\tilde{K}$ satisfying the same properties as $K$ satisfies that $\dim(E^s_{\tilde K}) \geq \dim(E^s_K)$ (this will be used only for Case C)).

\subsubsection*{Case A): There exists a chain recurrent central segment. $K$ has type $(R)$}

Assume that the set $K\en H$ admits a chain recurrent central segment. That is, there exists a curve $\gamma$ tangent to $E^c$ in a point of $K$, which is contained in $H$ and such that $\gamma$ is contained in a compact, invariant, chain transitive set in $U$, a small neighborhood of $K$.

In this case, the results of \cite{C1} (Addendum 3.9) imply that there are periodic orbits in the same chain recurrence class as $K$ (i.e. $H$) with index $\dim E^s \leq \alpha-1$, a contradiction.

\subsubsection*{Case B): $K$ has type $(N)$, $(H)$ or $(P_{SN})$}

 If $K$ has type $(H)$, one can apply Proposition 4.4 of \cite{C} which implies that there is a weak periodic orbit in $H$ giving option $2)$ of Theorem B.1.

Cases $(N)$ and $(P_{SN})$ give a family of central curves $\gamma_x$ $\forall x \in K$ (tangent to $E^c$, see \cite{C}) which satisfy that $f(\gamma_x)\en \gamma_{f(x)}$. It is not difficult to see that there is a neighborhood $U$ of $K$ such that for every invariant set in $U$ the same property will be satisfied (see remark 2.3 of \cite{C1}).

Consider a set $\umb{K}=K\cup \bigcup_n \mathcal O_n$ where $\mathcal O_n$ are close enough periodic orbits converging in Hausdorff topology to $K$ (these are given, for instance, by \cite{C3}) which we can suppose are contained in $U$.

So, since for some $x\in K$, we have that $W^{uu}_{loc}(x)$ will intersect $W^{cs}_{loc}(p_n)$ in a point $z$ (for a point $x$, the local center stable set, $W^{cs}_{loc}(x)$ is the union of the local strong stable leaves of the points in $\gamma_x$).

Since the $\omega$-limit set of $z$ must be a periodic point (see Lemma 3.13 of \cite{C1}) and since $H$ is Lyapunov stable we get that there is a periodic point of index $\alpha$ which is weak, or a periodic point of smaller index in $H$ which gives a contradiction.

\subsubsection*{ Case C): $K$ has type ($P_{UN}$) or ($P_{SU}$)}

One has a minimal set $K$ which is contained in a Lyapunov stable homoclinic class and it admits a partially hyperbolic splitting with one dimensional center with zero exponents and type $(P_{UN})$ or $(P_{SU})$.

This gives that given a compact neighborhood $U$ of $K$, there exists a family of $C^1$-curves $\gamma_x:[0,1]\to U$  ($\gamma_x(0)=x$) tangent to the central bundle such that $f^{-1}(\gamma_x([0,1]))\en \gamma_{f^{-1}(x)}([0,1))$. This implies that the preimages of these curves remain in $U$ for past iterates and with bounded length.

They also verify that the chain unstable set of $K$ restricted to $U$ (that is, the set of points that can be reached from $K$ by arbitrarily small pseudo orbits contained in $U$) contains these curves. Since $H$ is Lyapunov stable, this implies that these curves are contained in $H$.

Assume we could extend the partially hyperbolic splitting from $K$ to a dominated splitting $T_{K'}M=E_1\oplus E^c \oplus E_3$ in a chain transitive set $K'\en H$ containing $\gamma_x([0,t))$ for some $x\in K$ and for some $t\in (0,1)$.

Since the orbit of $\gamma_x([0,1])$ remains near $K$ for past iterates, we can assume (by choosing $U$ sufficiently small) that the bundle $E_3$ is uniformly expanded there. So, there are uniformly large unstable manifolds for every point in $\gamma_x([0,1])$ and are contained in $H$.

If we prove that $E_1$ is uniformly contracted in all $K'$, since we can approach $K'$ by weak periodic orbits, we get weak periodic in the class since its strong stable manifold (tangent to $E_1$) will intersect $H$.

To prove this, we use that for $K$ the dimension of $E^s$ is minimal. So $E_1$ must be stable, otherwise, we would get that, using Theorem B.1 again, there is a partially hyperbolic minimal set inside $K'$ with stable bundle of dimension smaller than the one of $K$, a contradiction.

The fact that we can extend the dominated splitting and approach the point $y$ in $\gamma_x((0,1))$ by weak periodic points is given by  Lemma B.1 below.

This concludes the proof of Theorem 5.1.
\lqqd

The following Lemma allows to conclude the proof of the previous Theorem but it can also have interest on its own (see for example \cite{C4} chapter 9). Its proof is postponed to the next subsection.

\begin{lemaB1}\label{encontrarpuntos} There exists a residual subset $\Res' \en \diff 1 \backslash \overline{\Tang}$ such that every $f\in \Res'$ verifies the following. Given a compact invariant set $K$ such that

\begin{itemize}
\item[-]$K$ is a chain transitive set.
\item[-]$K$ admits a partially hyperbolic splitting $T_KM = E^s \oplus E^c \oplus E^u$ where $E^s$ is uniformly contracting, $E^u$ is uniformly expanding and $\dim E^c= 1$.
 \item[-] Any invariant measure supported in $K$ has zero Lyapunov exponents along $E^c$.
\end{itemize}

Then, for every $\delta>0$, there exists $U$, a neighborhood of $K$ such that for every $y\in U$ satisfying:
\begin{itemize}
\item[-] $y$ belongs to the local chain unstable set $pW^u(K,U)$ of $K$ (that is, for every $\eps>0$ there exists an $\eps-$pseudo orbit from $K$ to $y$ contained in $U$)
    \item[-] $y$ belongs to the chain recurrence class of $K$
\end{itemize}

\noindent we have that there exist $p_n \to y$, periodic points, such that:

\begin{itemize}
\item[-] The orbit $\cO(p_n)$ of the periodic point $p_n$ has its $\dim E^s+1$ Lyapunov exponent contained in $(-\delta,\delta)$.
\item[-] For large enough $n_0$, if $\tilde K = K \cup \overline{\bigcup_{n>n_0} \mathcal O(p_n)}$, then we can extend the partially hyperbolic splitting to a dominated splitting of the form $T_{\tilde K}M =E_1 \oplus E^c \oplus E_3$.
\end{itemize}
\end{lemaB1}

\subsection*{B.1. Proof of Lemma B.1.}\label{Lema}

We shall first prove a perturbation result and afterwards we shall deduce Lemma B.1 with a standard Baire argument. One can compare this Lemma with Lemma 3.2 of \cite{Y} which is a slightly weaker version of this.

\begin{lemaB2}\label{perturbacion} There exists a residual subset $\Res \en \Diff^1(M)$ such that every $f\in \Res$ verifies the following. Given:
\begin{itemize}
\item[-] $K$ a compact chain transitive set.
\item[-] $U$ a neighborhood of $K$ and $y\in U$ verifying that $y$ is contained in the local chain unstable set $pW^u(K,U)$ of $K$ and in the chain recurrence class of $K$.
\item[-] $\U$ a $C^1$-neighborhood of $f$.
\end{itemize}
Then, there exists $l>0$ such that, for every $\nu>0$ and $L>0$ we have $g\in \U$ with a periodic orbit $\cO$ with the following properties:
\begin{itemize}
\item[-] There exists $p_1\in \cO$ such that $d(f^{-k}(y),g^{-k}(p_1))<\nu$ for every $0\leq k \leq L$.
\item[-] There exists $p_2\in \cO$ such that $\cO \backslash \{p_2, \ldots, g^l(p_2)\} \en U$.
\end{itemize}
\end{lemaB2}

\dem{\!\!} The argument is similar as the one in section 1.4 of \cite{C}. We must show that after an arbitrarily small perturbation, we can construct such periodic orbits.

Let $\Res$ be the $C^1$-residual subset given by Theorem A.1. Let $f\in \Res$.

Consider a point $y$ as above.  We can assume that $y$ is not chain recurrent in $U$, otherwise $y$ we would be accumulated by periodic orbits contained in $U$ (see property $a6)$ of Theorem A.1) and that would conclude without perturbing.

For every $\eps>0$ we consider an $\eps-$pseudo orbit $Y_\eps = (z_0 , z_1, \ldots, z_n)$ with $z_0 \in K$ and $z_n=y$ contained in $U$. Using that $y$ is not chain recurrent in $U$, we get that for $\eps$ small enough we have that $B_\nu(y)\cap Y_\eps = \{y\}$ where $B_\nu(y)$ is the ball of radius $\nu$ and $\nu$ is small. So if we consider a Hausdorff limit of the sequence $Y_{1/n}$ we get a compact set $Z^{-}$ for which $y$ is isolated and such that it is contained in the chain unstable set of $K$ restricted to $U$. Notice that $Z^-$ is backward invariant.

If we now consider the pair $(\Delta^{-}, y)$ where $\Delta^{-} = Z^{-}\backslash \{f^{-n}(y)\}_{n\geq 0}$ we get a pair as the one obtained in Lemma 1.11 of \cite{C} where $y$ plays the role of $x^{-}$. Notice that $\Delta^{-}$ is compact and invariant.

Now, we consider $\tilde U \en U$ a small neighborhood of $K\cup \Delta^-$ such that $y\notin \tilde U$. Take $x \in H \cap U^c$ where $H$ is the chain recurrence class of $K$.

Consider $X_\eps =(z_0, \ldots z_n)$ an $\eps-$pseudo orbit such that $z_0=x$ and $z_n\in K$. Take $z_j$ the last point of $X_\eps$ outside $\tilde U$. Since we chose $\tilde U$ small we have $z_j \in U\backslash \tilde U$. We call $\tilde X_\eps$ to $(z_j,\ldots z_n)$.

Consider $Z^+$ the Haussdorff limit of the sequence $\tilde X_{1/n}$ which will be a forward invariant compact set which intersects $U\backslash \tilde U$. Since $y$ is not chain recurrent in $U$ we have that $y\notin Z^+$.

We consider a point $x^+ \in Z^+ \cap U\backslash \tilde U$. This point satisfies that one can reach $K$ from $x^+$ by arbitrarily small pseudo orbits. We get that the future orbit of $x^+$ does not intersect the orbit of $y$.

Consider $\U$ a neighborhood of $f$ in the $C^1$ topology. Hayashi's connecting lemma (Theorem A.2) gives us $N>0$ and neighborhoods $W^+ \en \hat{W}^+$ of $x^+$ and $W^- \en \hat{W}^-$ of $y$ which we can consider arbitrarily small so, we can suppose that

\begin{itemize}
\item[-] All the iterates $f^i(\hat{W}^+)$ and $f^{j}(\hat{W}^-)$ for $0\leq i,j \leq N$ are pairwise disjoint.
\item[-] The iterates $f^i(\hat{W}^+)$ for $0\leq i \leq N$ are disjoint from the past orbit of $y$.
\end{itemize}

Since there are arbitrarily small pseudo orbits going from $y$ to $x^+$ contained in $H$ and $f$ is generic, property $a6)$ of Theorem A.1 gives us an orbit $\{x_0, \ldots, f^l(x_0) \}$ in a small neighborhood of $H$ and such that $x_0 \in W^-$ and $f^l(x_0)\in W^+$.

The same argument gives us an orbit $\{x_1, \ldots, f^k(x_1)\}$ contained in $U$ such that $x_1 \in W^+$ and $f^k(x_1) \in W^{-}$. In fact, we can choose it so that $d(f^{-i}(y),f^{k-i}(x_1))<\nu$ for $0\leq i \leq L$ (this can be done using uniform continuity of $f^{-1}, \ldots, f^{-L}$ and choosing $f^k(x_1)$ close enough to $y$).

Using Hayashi's connecting lemma (Theorem A.2) we can then create a periodic orbit $\mathcal O$ for a diffeomorphism $g \in \U$ which is contained in $\{x_0,\ldots, f^l(x_0)\} \cup \hat{W}^+ \cup \ldots \cup f^N(\hat{W^+}) \cup \{x_1, \ldots , f^k(x_1)\} \cup \hat{W}^- \cup \ldots \cup f^N(\hat{W}^-)$ (in the proof of Proposition 1.10 of \cite{C}) is explained how one can compose two perturbations in order to close the orbit).

Notice that from how we choose $\hat{W}^+$ and $\hat{W}^-$ and the orbit $\{x_1, \ldots, f^k(x_1)\}$ we get that the periodic orbit we create with the connecting lemma satisfies that $d(f^{-i}(y), g^{-i}(p_n)) < \nu$ for $0 \leq i \leq L$ and some $p_n$ in the orbit. This is because $\{x_1, \ldots, f^{k-1}(x_1)\}$ does not intersect $\hat{W}^- \cup \ldots \cup f^N(\hat{W}^-)$ and $\{f^{k-L}(x_1), \ldots f^k(x_1)\}$ does not intersect $\hat{W}^+ \cup \ldots \cup f^N(\hat{W^+})$ (in fact, this gives that the orbit of $p_n$ for $g$ contains $\{f^{k-L}(x_1), \ldots f^k(x_1)\}$).

Also, since $\hat{W}^+ \cup \ldots \cup f^N(\hat{W^+}) \cup \{x_1, \ldots, f^k(x_1)\}  \en U$ we get that, except for maybe $l$ consecutive iterates of one point in the resulting orbit, the rest of the orbit is contained in $U$. This concludes the proof.

\lqqd

\dem{of Lemma B.1}  Take $x \in M$, an let $m, e, t \in \N$. We consider $\U(m,e, t, x)$ the set of $C^1$ diffeomorphisms $g\in \Diff^1(M)$ with a periodic orbit $\cO$ satisfying:

\begin{itemize}
\item[-] $\cO$ is hyperbolic.
\item[-] Its $e+1$ Lyapunov exponent of $\cO$ is contained in $(-1/m,1/m)$.
 \item[-] There exists $p\in \cO$ such that $d(p,x) < 1/t$.
\end{itemize}

This set is clearly open in $\Diff^1(M)$. Let $\{x_s\}$ be a countable dense set of $M$. We define $\Res_{m,e, t, s} =  \U(m,e, t, x_s) \cup \Diff^1(M)\backslash \overline{\U(m,e,t , x_s)}$ which is open and dense by definition. Consider $\Res_1 = \bigcap_{m,e, t, s} \Res_{m,e, t, s}$ which is residual. Finally, taking $\Res$ as in Lemma B.2, we consider

 $$\Res'=(\Res \cap \Res_1) \backslash \overline{\Tang}$$

\smallskip

Consider $K$ compact chain transitive and with a partially hyperbolic splitting $T_KM = E^s \oplus E^c \oplus E^u$ with $\dim E^c=1$. We assume that any invariant measure supported in $K$ has Lyapunov exponent equal to zero.

Choose $\delta>0$ small enough. Since $f$ is far from tangencies, a Theorem of Wen \cite{W2} (see also Corollary 1.1 of \cite{C}) gives us that every periodic orbit having its $\dim E^s +1$  Lyapunov exponent in $(-\delta, \delta)$ admits a dominated splitting $E_1 \oplus E^c \oplus E_3$ with uniform strength (that is, if there is a set $\{O_n\}$ of periodic orbits with their $\dim E^s+1$ Lyapunov exponents in $(-\delta,\delta)$, then the dominated splitting extends to the closure).

We choose $U_1$, an open neighborhood of $K$ such that every invariant measure supported in $U_1$ has its Lyapunov exponents in $(-1/2m_0, 1/2m_0)$ where $1/m_0 <\delta$.

We can assume that $U_1$ verifies that there are $Df$ invariant cones $\cC^{uu}$ and $\cC^{cu}$ around $E^u$ and $E^c\oplus E^u$ respectively, defined in $U_1$. Similarly, there are in $U_1$, $Df^{-1}$ invariant cones $\cC^{ss}$,  $\cC^{cs}$ around $E^s$ and $E^s\oplus E^c$ respectively\footnote{If the reader is not familiar with the relationship between cones and dominated splitting there is a short survey on this in \cite{NewhouseConos}.}.

 We can assume, by choosing an adapted metric (see \cite{HPS} or \cite{GourmelonAdaptada}), that for every $v\in \cC^{ss}$ we have $\|Df v \|< \lambda \|v\|$ and for every $v\in \cC^{uu}$ we have $\| Df v \|>\lambda^{-1} \|v\|$ for some $\lambda<1$. There exists $\U_1$, a $C^1$-neighborhood of $f$ such that for every $g\in \U_1$ the properties above remain true.

Given $U$ neighborhood of $K$ such that $\overline U \en U_1$ we have that any $g$ invariant set contained in $U$ admits a partially hyperbolic splitting.

We now consider $y\in pW^u(K,U)$ which is contained in the chain recurrence class of $K$.

\begin{claim}\label{AfirmacionConos} Given $t$, for any $x_s$ with $d(x_s,y)<1/2t$ we get that  $f\in \U(m_0,\dim E^s, t, s) $.
\end{claim}

\dem{of Claim \ref{AfirmacionConos}} Since $f$ is in $\Res_1'$ it is enough to show that every neighborhood of $f$ intersects $\U(m_0,\dim E^s, t, s)$.

Choose a neighborhood $\U$ of $f$ and consider $\U_0 \en \U$ given by Frank's Lemma (\cite{frankslema}) such that we can perturb the derivative of some $g\in \U_0$ in a finite set of points less than $\xi$ and obtain a diffeomorphism in $\U$.

For $\U_0$, Lemma \ref{perturbacion} gives us a value of $l<0$ such that for any $L>0$ and there exists $g_L \in \U_0$ and a periodic orbit $\cO_L$ of $g_L$ such that there is a point $p_1\in \cO_L$ satisfying that $d(g^{-i}(p_1), f^{-i}(y))<1/2t$ ($0\leq i \leq L$) and a point $p_2\in \cO_L$ such that $\cO_L \backslash \{ p_2 , \ldots, g^l(p_2)\}$ is contained in $U$. We can assume that $\cO$ is hyperbolic.

We must perturb the derivative of $\cO_L$ less than $\xi$ in order to show that the $\dim E^s +1$ Lyapunov exponent is in $(-1/m_0,1/m_0)$.

Notice that if we choose $L$ large enough, we can assume that the angle of the cone $Dg^L(\cC^{\sigma}(g^{-L}(p_2))$ is arbitrarily small ($\sigma= uu, cu$). In the same way, we can assume that the angle of the cone $Dg^{-L}(\cC^{\tilde \sigma}(g^{L+l}(p_2)$ is arbitrarily small ($\tilde \sigma=cs, ss$ respectively).

Since $l$ is fixed, we get that for $p \in \cO_L \cap U^c$ (if there exists any, we can assume it is $p_2$), it is enough to perturb less than $\xi$ the derivative in order to get the cones $Dg^L(\cC^{\sigma}(g^{-L}(p_2))$ and $Dg^{-L-l}(\cC^{\tilde \sigma}(g^{L+l}(p_2)$ transversal (for $\sigma=uu, cu$ and $\tilde \sigma=cs, ss$ respectively). This allows us to have a well defined dominated splitting above $\cO_L$ (which may be of very small strength) which in turn allow us to define the $\dim E^s+1$ Lyapunov exponent. Since the orbit $\cO_L$ spends most of the time inside $U$, and any measure supported in $U$ has its center Lyapunov exponent in $(-1/2m_0,1/2m_0)$ we get the desired property.

\finobs

Taking $t=n\to \infty$ and using Wen's result, we get a sequence of periodic points $p_n \to y$ such that if $\cO(p_n)$ are their orbits, the set $\tilde K = K \cup \overline{ \bigcup_n \cO(p_n)}$ admits a dominated splitting $T_{\tilde K}M = E_1 \oplus E^c \oplus E_3$ extending the partially hyperbolic splitting.

This concludes the proof of Lemma B.1.

\lqqd

\end{document}